\documentclass[letterpaper, 10 pt, conference]{styles/ieeeconf_new}
\IEEEoverridecommandlockouts                             
\usepackage{etoolbox}
\usepackage{cuted} 

\AfterEndEnvironment{strip}{\leavevmode}
\usepackage{url}
\usepackage[dvipsnames]{xcolor}
\usepackage{wrapfig}
\newif\iffigs
\figstrue 

\usepackage{amsmath}
\usepackage{float}
\usepackage{setspace}
\usepackage{graphicx}
\usepackage{mathrsfs}
\usepackage{amssymb}
\usepackage{nicefrac}
\usepackage{algorithm}
\usepackage{stfloats}
\usepackage[noend]{algpseudocode}

\makeatletter
\newcommand\fs@spaceruled{\def\@fs@cfont{\bfseries}\let\@fs@capt\floatc@ruled
  \def\@fs@pre{\vspace{0.4\baselineskip}\hrule height.8pt depth0pt \kern2pt}%
  \def\@fs@post{\vspace{-0.4\baselineskip}\kern2pt\hrule\relax\vspace{-12pt}}%
  \def\@fs@mid{\kern2pt\hrule\kern2pt}%
  \let\@fs@iftopcapt\iftrue}
\makeatother

\usepackage{listings}
\lstnewenvironment{itemlisting}[1][]
 {%
  \mbox{}
  \vspace*{-\baselineskip}
  \lstset{
    xleftmargin=\leftmargin,
    linewidth=\linewidth,
    #1
  }%
 }
 {}
\usepackage{flushend}
\usepackage{multicol}

\usepackage{hyperref}

\newcommand{\scaleForAppendixMath}{0.88}

\usepackage{lipsum}

\newtheorem{thm}{Theorem}[section]

\newtheorem{lem}[thm]{Lemma}

\title{\LARGE \bf Symmetric Stair Preconditioning of Linear Systems\\for Parallel Trajectory Optimization}

\author{Xueyi Bu$^{1}$, Brian Plancher$^{2}$
\thanks{This material is based upon work supported by the National Science Foundation (under Award 2246022). Any opinions, findings, conclusions, or recommendations expressed in this material are those of the authors and do not necessarily reflect those of the funding organizations.}
\thanks{$^{1}$Xueyi Bu is with the School of Engineering and Applied Science, Columbia University, New York, NY. {\tt\footnotesize xb2133@columbia.edu}}%
\thanks{$^{2}$Brian Plancher is with Barnard College, Columbia University, New York, NY. {\tt\footnotesize bplancher@barnard.edu}}%
}%

\begin{document}
\maketitle
\thispagestyle{empty}
\pagestyle{empty}


\begin{abstract}
    There has been a growing interest in parallel strategies for solving trajectory optimization problems. One key step in many algorithmic approaches to trajectory optimization is the solution of moderately-large and sparse linear systems. Iterative methods are particularly well-suited for parallel solves of such systems. However, fast and stable convergence of iterative methods is reliant on the application of a high-quality preconditioner that reduces the spread and increase the clustering of the eigenvalues of the target matrix. To improve the performance of these approaches, we present a new parallel-friendly symmetric stair preconditioner. We prove that our preconditioner has advantageous theoretical properties when used in conjunction with iterative methods for trajectory optimization such as a more clustered eigenvalue spectrum. Numerical experiments with typical trajectory optimization problems reveal that as compared to the best alternative parallel preconditioner from the literature, our symmetric stair preconditioner provides up to a 34\% reduction in condition number and up to a 25\% reduction in the number of resulting linear system solver iterations.

\end{abstract}

\section{Introduction} \label{sec:intro}
Trajectory optimization algorithms~\cite{Betts01} are a powerful, and state-of-the-art set of tools for synthesizing dynamic motions for complex robots~\cite{Tassa12,Xi14,Posa16,Apgar18,Howell19,Kuindersma23Talk}. There has been historical interest in parallel strategies~\cite{Betts91} for solving trajectory optimization problems and several more recent efforts have shown that significant computational benefits are possible by exploiting parallelism in different stages of the algorithm on heterogeneous hardware platforms including multi-core CPUs~\cite{Kouzoupis16,Giftthaler17,Farshidian17,Astudillo22}, GPUs~\cite{Antony17,Plancher18,Plancher19a,Pan19,Plancher21,Plancher22GRiDGPUAcceleratedRigidBodyDynamicsAnalyticalGradientsb}, and FPGAs~\cite{Sacks18,Neuman21,Neuman23RoboShape,yang2023rbdcore}. This shift toward parallelism is growing increasingly important with the impending end of Moore's Law and the end of Dennard Scaling, which have led to a utilization wall that limits the performance a single CPU chip can deliver~\cite{Esmaeilzadeh11,Venkatesh10}.

One key step in many algorithmic approaches to trajectory optimization is the solution of moderately-large and sparse linear systems. These systems can be solved by direct factorization or through iterative fixed point approaches. While many state-of-the-art solvers leverage factorization-based approaches~\cite{stellato2020osqp,frison2020hpipm}, iterative methods, like the Preconditioned Conjugate Gradient (PCG) algorithm~\cite{eisenstat1981efficient}, are particularly well-suited for parallelism, as they are computationally dominated by matrix-vector products and vector reductions~\cite{saadIterativeMethodsSparse2003,plancherGPUAccelerationRealTime2022}. Furthermore, in the context of very-large-scale linear systems GPU implementations of PCG, and other iterative methods, have already been shown to outperform both factorization and iterative approaches on the CPU~\cite{schubigerGPUAccelerationADMM2020,helfenstein2012parallel}. However, fast and stable convergence of such iterative methods is reliant on the use of a high-quality symmetric preconditioner that reduces the spread and increase the clustering of the eigenvalues of the target matrix~\cite{saadIterativeMethodsSparse2003,Nocedal06,Shewchuk94,liNewApproximateFactorization2009}.  

To improve the performance of these approaches, we present a new parallel-friendly symmetric stair preconditioner. We prove that our preconditioner has advantageous theoretical properties when used in conjunction with iterative methods for trajectory optimization such as a more clustered eigenvalue spectrum than previous parallel preconditioners. Numerical experiments with typical trajectory optimization problems reveal that as compared to the best alternative parallel preconditioner from the literature, our symmetric stair preconditioner provides up to a 34\% reduction in condition number and up to a 25\% reduction in the number of PCG iterations needed for convergence.

\section{Background and Related Work} \label{sec:bgrelated}
\subsection{Direct Trajectory Optimization}
\label{sec:background_trajopt}
Trajectory optimization~\cite{Betts01}, also known as numerical optimal control, solves an (often) nonlinear optimization problem to compute a robot's optimal path through an environment as a series of states (x $\in \mathbb{R}^n$) and controls (u $\in \mathbb{R}^m$). These problems assume a discrete-time dynamical system,
\begin{equation}
x_{k+1} = f(x_k,u_k,h), \quad x_0 = x_s, \label{eq:dynamics}
\end{equation}
with a timestep $h$, and minimize an additive cost function,
\begin{equation}
    J(X,U) = \ell_f(x_N) + \sum^{N-1}_{k = 0} \ell(x_k,u_k). \label{eq:cost}
\end{equation}

One approach for solving such problems is through the use of \emph{direct methods}. At each iteration of a direct method, a quadratic approximation of the nonlinear problem defined by Equations~\ref{eq:dynamics} and ~\ref{eq:cost} is formed around a nominal trajectory ($X,U$), resulting in a quadratic program (QP). The resulting KKT system~\cite{karush39minima,kuhn1951nonlinear} is then solved to update the nominal trajectory. This process is repeated until convergence.

One approach to solving the KKT system is through two step process using the symmetric positive definite \emph{Schur Complement}, $S$.\footnote{While $S$ is actually symmetric negative definite, one can instead negate it and solve for $-\lambda^*$, enabling us to simplify our proofs by treating it as symmetric positive definite in the remainder of this work.} This process first solves for the optimal Lagrange multipliers, $\lambda^*$, and then recovers the state and control update, $\delta z^* = [\delta x^*, \delta u^*]$ as follows (where $c$ and $C$ are the constraint value and Jacobian, and $g$ and $G$ are the cost Jacobian and Hessian respectively):
\begin{equation} \label{eq:schur}
\begin{split}
    S &= -CG^{-1}C^T        \hspace{20pt}      
    \gamma = c - CG^{-1}g \\
    S\lambda^* &= \gamma  \hspace{54pt}     
    \delta z^* = G^{-1}(g - C^T\lambda^*) \\
\end{split}
\end{equation}

If we denote the dynamics Jacobians as $A_k = f_{x_k}$ and $B_k = f_{u_k}$, and cost Jacobians and Hessians as $q_k = J_{x_k}$, $r_k = J_{u_k}$, $Q_k = J_{x_k x_k}$, and $R_k = J_{u_k u_k}$, we can define: 
\begin{equation}
\begin{aligned}
\theta_k & =-A_k Q_k^{-1} A_k^T-B_k R_k^{-1} B_k^T-Q_{k+1}^{-1} \\
\phi_k & =A_k Q_k^{-1} \\
\zeta_k & =A_k Q_k^{-1} q_k+B_k R_k^{-1} r_k-Q_{k+1}^{-1} q_{k+1}.
\end{aligned}
\end{equation}
The Schur complement then takes the following block-tridiagonal form: 
\begin{equation}
\begin{aligned}
& S=\left(\begin{matrix}
-Q_0^{-1} & \phi_0^T & & &  \\
\phi_0 & \theta_0 & \phi_1^T & & \\
 & & \ddots &  & \\
& &  \phi_{N-2} & \theta_{N-2} & \phi_{N-1}^T \\
& &   &\phi_{N-1} & \theta_{N-1} 
\end{matrix}\right) \\
&\gamma=c+\left(
-Q_0^{-1} q_0 \quad \zeta_0\quad \zeta_1 \quad \ldots \quad \zeta_{N-1}
\right)^T.
\end{aligned}
\end{equation}

For most trajectory optimization problems $Q$ and $R$ are chosen to be symmetric positive definite, or made symmetric positive definite by construction (e.g., in the case of quadratic penalty methods), and $C^T$ has full column rank. Furthermore, by construction $\theta$ is both invertible and symmetric. 
Similarly, although $A_k$ and $B_k$ will vary depending on the integrator in use, for most popular explicit and semi-implicit integrators they take the form
\begin{equation}
A = I + h M \quad B= h N,
\end{equation}
where $M$ and $N$ are functions of the system dynamics $f$. Additional conditions, such as Lipschitz continuity, are usually imposed on $f$ to ensure the local existence and uniqueness of the solution to the dynamical system~\cite{khalilNonlinearSystems2002}. As a result, $\Vert M \Vert$ is bounded; thus, for sufficiently small $h$, $h\Vert M \Vert < 1$ and $A = I + hM$, and thus $\phi$, is invertible. 
Finally, given all of the above, the dominant computational step in this approach is solving for $\lambda^*$ in Equation~\ref{eq:schur}.

\subsection{Parallel Iterative Solvers}
There has been a significant amount of prior work developing general purpose parallel iterative solvers, many of which target the GPU~\cite{saadIterativeMethodsSparse2003,Bolz03,Liu13,Anzt17,Anzt18,Flegar19,Schubiger20}. 
%
Iterative methods solve the problem $S\lambda^* = \gamma$ for a given $S$ and $\gamma$ by iterative refining an estimate for $\lambda$ up to some tolerance $\epsilon$. The most popular of these methods is the Conjugate Gradient (CG) method, which is applicable to systems where $S$ is positive definite. The convergence rate of CG is directly related to the spread of the Eigenvalues of $S \in \mathbb{R}^{n \times n}$, converging faster when they are clustered. We also note that clustered, and moderate in magnitude, eigenvalues also have the added benefit of avoiding round-off and overflow errors caused by iterative floating point math.

To improve the performance of CG, a symmetric preconditioning matrix $\Phi \approx S$ is often applied to instead solve the equivalent problem $\Phi^{-1} S\lambda^* = \Phi^{-1} \gamma$~\cite{pearson2020preconditioners}. The resulting preconditioned CG (PCG) algorithm leverages matrix-vector products with $S$ and $\Phi^{-1}$, as well as vector reductions, both parallel friendly operations. 
Finally, we note that PCG requires a symmetric preconditioner~\cite{pearson2020preconditioners}.


\subsection{Parallel Preconditioners}

There is a large literature on preconditioning~\cite{saadIterativeMethodsSparse2003,adamsStepPreconditionedConjugate1985, adamsAdditivePolynomialPreconditioners1989, concusNumericalSolutionNonlinear1978, muzhinjiOptimalBlockPreconditioner2023, benziPreconditioningTechniquesSaddle2008} and there are many different preconditioners that are optimized for computation on vector or parallel processors. 
The most popular of these are the Jacobi and Block-Jacobi preconditioners. These set:
\begin{equation} \label{eq:jacobi}
\Phi = \text{diag}(S) \hspace{20pt} \text{or} \hspace{20pt} \Phi = \text{block-diag}(S). 
\end{equation}
Previous GPU based iterative solvers mainly leveraged these preconditioners~\cite{Flegar19,Schubiger20}. For block banded matrices, alternating and overlapping block preconditioners have also been used in previous work. These methods compute $\Phi^{-1}$ as a sum of the inverse of block-diagonal matrices that compose $S$~\cite{saadIterativeMethodsSparse2003,Galligani94}. Finally, Polynomial splitting preconditioners have also found widespread usage~\cite{saadIterativeMethodsSparse2003,liNewApproximateFactorization2009}. These follow the pattern $S = \Psi - P$ and compute a preconditioner where:
\begin{equation} \label{eq:polynomial_precon}
    S^{-1} \approx \Phi^{-1} = (I + \Psi^{-1}P + (\Psi^{-1}P)^2 \dots)\Psi^{-1}.
\end{equation}
We note that these preconditioners are only valid where:
\begin{equation}
    \rho(\Psi^{-1}P ) < 1,
\end{equation}
which is called a convergent splitting and will guarantee convergence when used with the CG algorithm~\cite{saadIterativeMethodsSparse2003,Nocedal06,Shewchuk94,liNewApproximateFactorization2009}. 
Also, increasing the degree of the polynomial computes a better approximation of $S$ and improves the convergence rate of the resulting PCG algorithm. However, this requires more computation to compute the preconditioner and also often creates a preconditioner with a larger bandwidth, requiring more memory. For block-banded matrices, like $S$ in our problem, the values in the true inverse decay exponentially as one moves away from the diagonal~\cite{Demko84}, this creates a tradeoff between the accuracy and both the memory and computational complexity of the preconditioner. 

\section{Stair Preconditioners} \label{sec:stair}
One polynomial splitting for (symmetric) block tridiagonal matrices that has been shown to outperform Jacobi and Block-Jacobi methods is the stair based splitting~\cite{Li09,Li11}. This splitting comes in two types, left (type 1) and right (type 2), depending upon the direction of the stair. For a 3x3 symmetric block tridiagonal $S$:
\begingroup
\setlength\arraycolsep{3pt}
\begin{equation} \label{eq:schur3x3}
\begin{split}
    S = &\begin{bmatrix}
        D_1 & O_1 & 0 \\
        O_1^T & \;\; D_2 \;\; & O_2 \\
        0 & O_2^T & D_3 \\
    \end{bmatrix}
\end{split}
\end{equation}
The stair splittings are:
\begingroup
\setlength\arraycolsep{3pt}
\begin{equation} \label{eq:stair}
\begin{split}
    \Psi_l = &\begin{bmatrix}
        {\color{Blue} D_1} & {\color{OliveGreen} 0} & 0 \\
        {\color{Blue} O_1^T} & \;\; {\color{Blue} D_2} \;\; & {\color{Blue} O_2} \\
        0 & {\color{OliveGreen} 0}& {\color{Blue} D_3} \\
    \end{bmatrix} \quad
    P_l = -\begin{bmatrix}
        0 & {\color{OliveGreen} O_1} & 0\\
        \;\; 0 \;\;  & 0 & \; 0 \; \\
        0 & {\color{OliveGreen} O_2^T} & 0 \\
    \end{bmatrix}\\
    \Psi_r = &\begin{bmatrix}
        {\color{Blue} D_1} & {\color{Blue} O_1} & 0 \\
        {\color{OliveGreen} 0} & \;\; {\color{Blue} D_2} \;\; & {\color{OliveGreen} 0} \\
        0 & {\color{Blue} O_2^T} & {\color{Blue} D_3} \\
    \end{bmatrix} \quad
    P_r = -\begin{bmatrix}
        0 & 0 & 0\\
        {\color{OliveGreen} O_1^T} & \;\; 0 \;\; & {\color{OliveGreen} O_2} \\
        0 & 0 & 0 \\
    \end{bmatrix}.\\
\end{split}
\end{equation}
\endgroup
The stair is invertible if and only if its diagonal blocks, $D$, are invertible. If so, its inverse has the same stair shaped sparsity pattern drawing from the same or neighboring block-rows of $S$, enabling efficient parallel computation~\cite{plancherGPUAccelerationRealTime2022,Lu99StairMatricesTheirGeneralizationsApplicationsIterativeMethodsGeneralizationSuccessiveOverrelaxationMethod}:
\begingroup
\setlength\arraycolsep{3pt}
\begin{equation} \label{eq:stair_inv}
\begin{split}
    \Psi^{-1} &= D^{-1}(2D - \Psi)D^{-1}\\
    \Psi_l^{-1} &= \begin{bmatrix}
        {\color{Blue} D_1^{-1}} & {\color{OliveGreen} 0} & 0 \\
        {\color{Blue} -D_2^{-1} O_1^T D_1^{-1}} & \;\; {\color{Blue} D_2^{-1}} \;\; & {\color{Blue} -D_2^{-1} O_2 D_3^{-1}} \\
        0 & {\color{OliveGreen} 0} & {\color{Blue} D_3^{-1}} \\
    \end{bmatrix}\\
    \Psi_r^{-1} &= \begin{bmatrix}
        {\color{Blue} D_1^{-1}} & {\color{Blue} -D_1^{-1} O_1 D_2^{-1}} & 0 \\
        {\color{OliveGreen} 0} & \;\; {\color{Blue} D_2^{-1}} \;\; & {\color{OliveGreen} 0} \\
        0 & {\color{Blue} -D_3^{-1} O_2^T D_2^{-1}} & {\color{Blue} D_3^{-1}} \\
    \end{bmatrix}.
\end{split}
\end{equation}
\endgroup
While, the left and right stair preconditioners are not symmetric, and thus cannot be used in the context of PCG, the \textit{additive polynomial preconditioner} can instead be used:
\begin{equation}
    \Phi_{add}^{-1} = \frac12 (\Psi_l^{-1} +\Psi_r^{-1}).
\end{equation}
However, as we shall see later on, this will result in a declustered spectrum that could negatively impact the performance of PCG.
In the remainder of this section we prove a few lemmas and theorems about the eigenpairs of $\Psi_l^{-1}$, $\Psi_r^{-1}$, and $\Phi_{add}^{-1}S$ which we will need for later proofs.\footnote{We note that throughout this paper we analyze the spectrum of $\Psi^{-1}S$ even though the resulting matrix is not guaranteed to be symmetric positive definite. This is possible because through the Cholesky factorization, $\Psi = L L^T$, we can instead from the symmetric positive definite $L^{-1} S L^{-T}$, and the matrices $L^{-1} S L^{-T}$ and $\Psi^{-1}S$ are similar.}



\subsection{Eigenpairs of $\Psi_l$ and $\Psi_l$}

We first prove a few lemmas about the eigenpairs of the left and right stair preconditioners.\footnote{For more details on stair preconditioners please see~\cite{Li09} and~\cite{Li11}.}

\medskip
\begin{lem} \label{lemma}
Given the stair-splittings of a symmetric block tridiagonal matrix $S=\Psi_l - P_l = \Psi_r - P_r$ for a $n\times n$ block $S$, where each block is $m\times m$, 
$\Psi_l^{-1}P_l$ and $\Psi_r^{-1}P_r$ have the same spectrum.
\end{lem}
\medskip

\begin{proof}
%
If $(\lambda, v)$ is an eigenpair of $P_r \Psi_r^{-1}$ and $\lambda \neq 0$:
\begin{equation}
    \Psi_r^{-1} P_r (\Psi_r^{-1} v) = \Psi_r^{-1} (\lambda v),
\end{equation}
and $(\lambda, \Psi_r^{-1}v)$ is an eigenpair of $\Psi_r^{-1}P_r$. 
Since $P_r$ is not invertible, $0$ is an eigenvalue for both $P_r \Psi_r^{-1}$ and $\Psi_r^{-1} P_r$.
Also as $P_r$ and $\Psi_r^{-1}$ are equal dimensioned square matrices then by Theorem 1.3.22 of~\cite{horn2012matrix}, $P_r\Psi_r^{-1}$ and $\Psi_r^{-1} P_r$ have the same spectrum.
Finally, $(\Psi_l^{-1}P_l)^T = P_r \Psi_r^{-1}$ and every matrix has the same spectrum as its transpose.
\end{proof}

\medskip
\begin{lem} \label{lemma}
For $v^T = (v_1^T, v_2^T, \dots, v_i^T, \dots, v_n^T) \in \mathbb{R}^{nm}$ and $v_i \in \mathbb{R}^{m}$ for $i = 1,2,\dots,n$,  we denote $v_e^T  = (0, v_2^T,\dots, 0,  v_{2j}^T, \dots)$ and $v_o^T = (v_1^T,0, \dots, v_{2j+1}^T, 0, \dots)$ such that $v = v_e + v_o$.

If $(\lambda, v = v_e + v_o)$ is an eigenpair of $\Psi_l^{-1}P_l$ and $\lambda \neq 0$ , then $(\lambda, u = v_e + \lambda v_o)$ is an eigenpair of $\Psi_r^{-1}P_r$.

If $(\lambda, u = u_e + u_o)$ is an eigenpair of $\Psi_r^{-1}P_r$ and $\lambda \neq 0$, then $(\lambda, v = \lambda u_e + u_o)$ is an eigenpair of $\Psi_l^{-1}P_l$.
\end{lem}
\medskip

\begin{proof}
$\Psi_l^{-1}P_l v_o = 0$ for all $v_o$ because of the zero columns of $\Psi_l^{-1}P_l$ (see Equation~\ref{eq:lz} in the appendix). Similarly, $P_r v_e = 0$ and $\Psi_r^{-1}P_r v_e = 0$ (see Equation~\ref{eq:rz} in the appendix).
\medskip

Suppose $(\lambda, v = v_e + v_o)$ is an eigenpair of $\Psi_l^{-1}P_l$, then 
\begin{equation} \label{eq:lwhole}
  \Psi_l^{-1}P_l v_e = \Psi_l^{-1}P_l v  = \lambda (v_e + v_o).
\end{equation}

If $\lambda \neq 0$, we split Equation~\ref{eq:lwhole} into two parts,
\begin{equation} \label{eq:split_1}
 (\Psi_l^{-1} - D^{-1})P_l v = (\Psi_l^{-1} - D^{-1})P_l v_e = \lambda v_e
\end{equation}
\begin{equation} \label{eq:split_2}
 D^{-1}P_l v =  D^{-1}P_l v_e = \lambda v_o \quad \;\;
\end{equation}
where $D$ is the diagonal blocks of $S$ (see Equations~\ref{eq:(l-d)p} and \ref{eq:dpl} in the appendix for more details).

Since $P_r = \Psi_r - S = D - \Psi_l$, then:
\begin{equation}
\begin{split}
D^{-1}P_r  &= I - D^{-1}\Psi_l =  (\Psi_l^{-1} - D^{-1})\Psi_l.
\end{split}
\end{equation}

This means that:
\begin{equation}
\begin{split}
D^{-1}P_r v &= (\Psi_l^{-1} - D^{-1})\Psi_l v \\
& = (\Psi_l^{-1} - D^{-1})\Psi_l (\frac1\lambda \Psi_l^{-1}P_l v) \\
& = \frac1\lambda (\Psi_l^{-1} - D^{-1}) P_l v \\
& = v_e \label{eq:o}
\end{split}
\end{equation}

Using the inverse of the stair matrix (Equation~\ref{eq:stair_inv}):
\begin{equation}
    \Psi_r^{-1} - D^{-1} = D^{-1} (D - \Psi_r)D^{-1} = D^{-1} P_lD^{-1}
\end{equation}
\begin{equation}
\begin{split}
    (\Psi_r^{-1} - D^{-1}) P_r v &= D^{-1} P_l D^{-1}P_r v \\
    &= D^{-1} P_l v_e \\
    &= \lambda v_o \label{eq:e}
\end{split}
\end{equation}

Combining Equations~\ref{eq:o} and~\ref{eq:e} as well as make use of the fact $P_r v_e = 0$, we obtain the following result:
\begin{equation}
\begin{split}
\Psi_r^{-1}P_r ( v_e + \lambda v_o) & = \Psi_r^{-1}P_r (\lambda  v) \\
    & = \lambda (D^{-1}P_r v + (\Psi_r^{-1} - D^{-1}) P_r  v) \\
    & = \lambda ( v_e + \lambda v_o) \label{eq:lm1}
\end{split}
\end{equation}

The second statement follows from Lemma 4.1 or can be proved with the same strategy as above.
\end{proof}
\medskip

Finally, we note that the eigenvectors with eigenvalues at $0$ can be divided into two types. The first type of these are the eigenvectors $v_o$ and $u_e$ which are formed because of the zero columns of $\Psi_l^{-1}P_l$ and $\Psi_r^{-1}P_r$ (see Equations~\ref{eq:lz} and~\ref{eq:rz} in the appendix). The second type of these are the eigenvectors $v_e$ and $u_o$ such that $P_l v_e = 0$ and $P_r u_o = 0$.

If $n$ is even, then $\Psi_l^{-1}P_l$ and $\Psi_r^{-1}P_r$ have exactly $m \frac n2 $ eigenpairs of the first type. If $n$ is odd, $\Psi_l^{-1}P_l$ will have $m \frac{n+1}{2}$ while $\Psi_r^{-1}P_r$ will have $m$ less such eigenpair and, instead, have $m$ more eigenpairs in the form of $(0, u_o)$ since, in this case, $P_r u_o = 0$, which has non-zero solutions even when all $O$s are invertible:
\begin{equation}
\left\{\begin{array}{lr}
        O_1^T u_1 + O_2 u_3 = 0\\
        O_3^T u_3 + O_4 u_5 = 0\\
        \cdots \\
         O_{n-2}^T u_{n-2} + O_{n-1} u_n = 0
        \end{array}\right.
\end{equation}

In the case of most integrators, as discussed before, $O$s are invertible, and $P_l v_e = 0$ if and only if $v_e = 0$. Hence, there are exactly $m \lceil\frac n2\rceil$ eigenpairs with eigenvalues at $0$ and $m \lfloor\frac n2\rfloor$ eigenpairs with non-zero eigenvalues. 


\subsection{Eigenpairs of $\Phi_{add}$}
We next construct the $mn$ eigenpairs of $\Phi_{add}^{-1}S$ from those of $\Psi_l^{-1}P_l$ and $\Psi_r^{-1}P_r$.
\medskip

\begin{thm}\label{thm}
If $(0, v_e)$ is an eigenpair of $\Psi_l^{-1}P_l$ then $(1, v_e)$ is an eigenpair of $\Phi_{add}^{-1}S$ and if $(0, u_o)$ is an eigenpair of $\Psi_r^{-1}P_r$ then $(1, u_o)$ is an eigenpair of $\Phi_{add}^{-1}S$.
\end{thm}
\medskip

\begin{proof}
Let $(0, v_e)$ be an eigenpair of $\Psi_l^{-1}P_l$, then:
\begin{equation}
\Psi_l^{-1} S v_e = \Psi_l^{-1} (\Psi_l - P_l)v_e = v_e - 0 v_e = v_e
\end{equation}
\begin{equation} \label{eq:r0}
\Psi_r^{-1} S v_e = \Psi_r^{-1} (\Psi_r - P_r)v_e = v_e - 0 = v_e.
\end{equation}
\begin{equation}
    \Phi_{add}^{-1}S v_e = \frac12 (\Psi_l^{-1} +\Psi_r^{-1})S v_e =\frac{1}{2} (v_e +v_e) = v_e.
\end{equation}

Let $(0, u_o)$ be an eigenpair of $\Psi_r^{-1}P_r$, then
\begin{equation}
\Psi_l^{-1} S u_o = \Psi_l^{-1} (\Psi_l - P_l)u_o = u_o - 0 = u_o
\end{equation}
\begin{equation}
\Psi_r^{-1} S u_o = \Psi_r^{-1} (\Psi_r - P_r)u_o = u_o - 0 u_o = u_o.
\end{equation}
\begin{equation}
    \Phi_{add}^{-1}S u_o = u_o.
\end{equation}
\end{proof}
\medskip

\begin{thm}\label{thm3.4}
Let $(\lambda, v = v_e + v_o)$ be an eigenpair of $\Psi_l^{-1}P_l$.
If $\lambda \neq 0$, then $(1-\frac12(\lambda\pm\sqrt{\lambda}),  v_e \pm \sqrt{\lambda} v_o)$ are eigenpairs of $\Phi_{add}^{-1}S$.
\end{thm}
\medskip

\begin{proof}
Let $(\lambda, v = v_e + v_o)$ be an eigenpair of $\Psi_l^{-1}P_l$: 
Since $\Psi_l^{-1}P_l v_o = 0$, 
\begin{equation} \label{eq:lvo}
    \Psi_l^{-1}P_l v_e = \lambda (v_e + v_o).
\end{equation}

Similarly, Equation~\ref{eq:lm1} and $\Psi_r^{-1}P_r v_e = 0$ suggests: 
\begin{equation} \label{eq:rve}
    \Psi_r^{-1}P_r v_o = v_e + \lambda v_o.
\end{equation}

We can therefore make an educated guess that the eigenvalues of $\Phi_{add}^{-1}S = \frac12 (\Psi_l^{-1} +\Psi_r^{-1})S =\frac{1}{2} (\Psi_l^{-1}S + \Psi_r^{-1}S)$ and should be of the form $u = a v_e + b v_o$. 
Following from Equations~\ref{eq:lvo} and~\ref{eq:rve}, we can therefore say that: 
\begin{equation} \label{eq:phi_add_S_u}
\scalebox{0.99}{$
\begin{split}
    \Phi_{add}^{-1}S u &= \frac{1}{2} (\Psi_l^{-1}S (a v_e + b v_o) + \Psi_r^{-1}S (a v_e + b v_o)) \\
    &= u - \frac{1}{2} \left[b\Psi_l^{-1}P_l v_e  + a\Psi_r^{-1}P_r v_o\right] \\
    &= u - \frac{1}{2} \left[b\lambda(v_e + v_o)  + a( v_e+ \lambda v_o)\right]\\
    &= u - \frac12 \left[(a+\lambda b) v_e + \lambda (a +  b) v_o\right] \\
\end{split}
$}
\end{equation}

Given our assumed form $u = a v_e + b v_o$, Equation~\ref{eq:phi_add_S_u} should simplify to: $\beta a v_e + \beta b v_o$. Thus, the eigenpair $(\lambda, u)$ must satisfy the linear system:
\begin{equation}
\frac12
    \left (\begin{matrix}
        1 & \lambda \\
        \lambda & \lambda 
    \end{matrix}
    \right)
    \left (\begin{matrix}
       a \\
        b
    \end{matrix}
    \right) = (1-\beta) \left (\begin{matrix}
       a \\
        b
    \end{matrix}
    \right)
\end{equation}

This means we can conclude that:
\begin{equation} \label{eq:lambda_beta}
\beta = 1 - \frac12(\lambda \pm \sqrt{\lambda}), \quad \quad b = \pm \sqrt{\lambda} a
\end{equation}
\end{proof}
\medskip

\begin{thm}\label{thm3.5}
If $S$ is symmetric block tridiagonal and positive definite, then $\Phi_{add}^{-1}S$ has real positive eigenvalues.
\end{thm}
\medskip

\begin{proof}
%
If $S$ is positive definite, then its diagonal block matrix, $D$, is positive definite and thus:
\begin{equation}
\begin{split}
    v^T S v & > 0 \\
    v^T (\Psi_l + \Psi_r - D) v &>0 \\
    v^T (\Psi_l + \Psi_r) v &>0
\end{split}
\end{equation}
    
Since $\Psi_l = \Psi_r^T$ we can also note that $\Psi_l$, $\Psi_r$, $\Psi_l^{-1}$, and $\Psi_r^{-1}$ are all positive definite as:
\begin{equation}
    v^T (\Psi_l + \Psi_r) v = 2v^T \Psi_l v =2v^T \Psi_r v
\end{equation}

As $(\Psi_l^{-1})^T = \Psi_r^{-1}$, $\Phi_{add}^{-1} = \frac12(\Psi_l^{-1} + \Psi_r^{-1})$ is symmetric and positive definite. 
$\Phi_{add}^{-1}S$ is then similar to the symmetric positive definite matrix $\frac12 (\Psi_l^{-1} + \Psi_r^{-1})^{1/2} S (\Psi_l^{-1} + \Psi_r^{-1})^{1/2}$ and must also have positive real eigenvalues.
We note that this is a special case of Theorem 3.4 of \cite{liNewApproximateFactorization2009} and Theorem 2.2 and Corollary 2.3 of \cite{adamsAdditivePolynomialPreconditioners1989}.
\end{proof}

\section{The Symmetric Stair Precontioner} \label{sec:symstair}
In this section, we introduce a new symmetric parallel preconditioner for symmetric positive definite block tridiagonal matrices which improves upon existing stair based preconditoners, the \emph{symmetric stair preconditioner}: 
\begin{equation}
    \Phi^{-1}_{sym} = \Psi_l^{-1} + \Psi_r^{-1} - D^{-1}.
\end{equation}
This new preconditioner can also simply be thought of as taking $\Psi_l^{-1}$ or $\Psi_r^{-1}$ and copying the off diagonal blocks across the diagonal and takes the form:
\begingroup
\setlength\arraycolsep{-1pt}
\begin{equation}
\scalebox{0.96}{$
    \Phi^{-1}_{sym} = \begin{bmatrix}
        {\color{Blue} D_1^{-1}} & {\color{Blue} -D_1^{-1} O_1 D_2^{-1}} & 0 \\
        {\color{Blue} -D_2^{-1} O_1^T D_1^{-1}} & \;\; {\color{Blue} D_2^{-1}} \;\; & {\color{Blue} -D_2^{-1} O_2 D_3^{-1}} \\
        0 & {\color{Blue} -D_3^{-1} O_2^T D_2^{-1}} & {\color{Blue} D_3^{-1}} \\
    \end{bmatrix}\\
$}
\end{equation}
\endgroup
In the remainder of this section, we will prove that $\Phi_{sym}^{-1}S$ is admissible for use with PCG, and that $\Phi_{sym}^{-1}S$ has a more clustered spectrum than $\Phi_{add}^{-1}S$.

We begin by showing that $\Phi_{sym}^{-1}S$ and $\Psi_l^{-1}S$ (or $\Psi_r^{-1}S$) share the same eigenvalues. And, that the multiplicities of the eigenvalues of $\Phi_{sym}^{-1}S$ that are less than one are doubled as compared to $\Psi_l^{-1}S$ (or $\Psi_r^{-1}S$).

\medskip
\begin{thm}\label{thm4.1}
If $(0, v_e)$ is an eigenpair of $\Psi_l^{-1}P_l$ then $(1, v_e)$ is an eigenpair of $\Phi_{sym}^{-1}S$. 
If $(0, u_o)$ is an eigenpair of $\Psi_r^{-1}P_r$ then $(1, u_o)$ are eigenpairs of $\Phi_{sym}^{-1}S$.
\end{thm}
\medskip

\begin{proof}
Let $(0, v_e)$ be an eigenpair of $\Psi_l^{-1}P_l$, then by Equations~\ref{eq:split_1},~\ref{eq:r0} and~\ref{eq:dpl} in the appendix:
\begin{equation}
\scalebox{0.88}{$
\begin{split}
    \Phi_{sym}^{-1}S v_e &=  (\Psi_r^{-1} + (\Psi_l^{-1}- D^{-1}))S v_e\\
    &= v_e + (\Psi_l^{-1}- D^{-1})(\Psi_l - P_l) v_e\\
    &= v_e  + (\Psi_l^{-1}- D^{-1})\Psi_l v_e - (\Psi_l^{-1}- D^{-1})P_l v_e\\
    &= v_e  + 0 - 0 v_e = v_e
\end{split}
$}
\end{equation}

Similarly, let $(0, u_o)$ be an eigenpair of $\Psi_r^{-1}P_r$, then:
\begin{equation}
\scalebox{0.88}{$
\begin{split}
    \Phi_{sym}^{-1}S u_o &=  (\Psi_l^{-1} + (\Psi_r^{-1}- D^{-1}))S u_o\\
    &= u_o + (\Psi_r^{-1}- D^{-1})(\Psi_r - P_r) u_o\\
    &= u_o + (\Psi_l^{-1}- D^{-1})\Psi_r u_o - (\Psi_r^{-1}- D^{-1})P_r u_o \\
    &= u_o  + 0 - 0 u_o = u_o
\end{split}
$}
\end{equation}
\end{proof}
\medskip

\medskip
\begin{thm}\label{thm4.2}
If $(\lambda, v = v_e + v_o)$ is an eigenpair of $\Psi_l^{-1}P_l$,
then $(1-\lambda, v_o)$ and $(1-\lambda, v_e)$ are eigenpairs of $\Phi_{sym}^{-1}S$.
\end{thm}
\medskip

\begin{proof}
Using $P_r v_e=0$, Equations~\ref{eq:split_1} and~\ref{eq:dpl}:
\begin{equation}
\scalebox{0.88}{$
\begin{split}
    \Phi_{sym}^{-1}S v_e &=  (\Psi_r^{-1} + (\Psi_l^{-1}- D^{-1}))S v_e\\
    &= \Psi_r^{-1}(\Psi_r - P_r) v_e + (\Psi_l^{-1}- D^{-1})(\Psi_l - P_l) v_e\\
    &= v_e + (\Psi_l^{-1}- D^{-1})\Psi_l v_e  - (\Psi_l^{-1}- D^{-1})P_l v_e \\
    &= v_e  + 0 - \lambda v_e\\
    &= (1-\lambda)v_e.
\end{split}
$}
\end{equation}

Hence, $(1-\lambda, v_e)$ is an eigenpair of $\Phi_{sym}^{-1}S$. The proof for $(1-\lambda, v_o)$ 
follows the same structure.
\end{proof}
\medskip

Finally, we will prove that $\Phi_{sym}^{-1}S$ is admissible for PCG, and that $\Phi_{sym}^{-1}S$ has a smaller condition number than $\Phi_{add}^{-1}S$.

\medskip
\begin{thm}\label{thm4.3}
If $S$ is symmetric block tridiagonal and positive definite, then:
\begin{enumerate}
    \item $\Phi_{sym}^{-1}S$ has real eigenvalues.
    \item $\Phi_{sym}^{-1} P$ has non-negative real eigenvalues and $\rho(\Phi_{sym}^{-1} P) < 1$.
    \item $\Phi_{sym}^{-1} S$ is positive definite and $\lambda(\Phi_{sym}^{-1} S) \in (0, 1]$.
\end{enumerate}
\end{thm}
\medskip

\begin{table*}[!b]
\begin{equation} \label{eq:max}
    \lambda_{\max}(\Phi_{sym}^{-1}S) = 1 - \lambda_{\min}(\Psi^{-1}P) \leq 1 < 1 - \frac12\left(\lambda^{+}(\Psi^{-1}P) - \sqrt{\lambda^{+}(\Psi^{-1}P)}\right) \leq \lambda_{\max}(\Phi_{add}^{-1}S) \leq \frac98
\end{equation}
\begin{equation} \label{eq:min}
\lambda_{\min}(\Phi_{sym}^{-1}S) = 1 - \lambda^{+}_{\max}(\Psi^{-1}P)> 1 - \frac12\left(\lambda^{+}_{\max}(\Psi^{-1}P) + \sqrt{\lambda^{+}_{\max}(\Psi^{-1}P)}\right) = \lambda_{\min}(\Phi_{add}^{-1}S) > 0
\end{equation}
\end{table*}

\begin{proof}
We prove each part in order:
\begin{enumerate}
    \item As $S$ is positive definite, then $\Phi_{sym}^{-1}S$ is similar to the symmetric matrix $S^{1/2}\Phi_{sym}^{-1}S^{1/2}$. Therefore, the eigenvalues of $\Phi_{sym}^{-1}S$ are real.
    \item As $\lambda \in \mathbb{R}$ from point one, then as we again derive $\lambda$ from $\Psi$, by Equation~\ref{eq:lambda_beta} and \hyperref[thm3.5]{Theorem 3.5}, $1-\frac12(\lambda\pm\sqrt{\lambda})$ will be real and positive if and only if $0 \leq \lambda < 1$. Thus, $\Phi_{sym}^{-1} P$ has non-negative real eigenvalues and $\rho(\Phi_{sym}^{-1} P) < 1$. This also shows that the stair splitting is convergent for a symmetric positive definite $S$.
    \item From points one and two and \hyperref[thm4.1]{Theorem 4.1} and \hyperref[thm4.2]{Theorem 4.2}, $1-\lambda \in (0, 1]$ and $\Phi_{sym}^{-1} S$ is positive definite.
    \end{enumerate}
\end{proof}

As a result of \hyperref[thm3.4]{Theorem 3.4}, \hyperref[thm3.5]{Theorem 3.5}, and \hyperref[thm4.3]{Theorem 4.3}, we have the following relationships between the maximum (Equation~\ref{eq:max}) and minimum (Equation~\ref{eq:min}) eigenvalues of the two preconditioners,
where $\lambda^{+}_{\max}(\Psi_x^{-1}P)$ and $\lambda^{+}_{\min}(\Psi_x^{-1}P)$ denotes the largest and smallest non-zero eigenvalues of $\Psi_x$. As we can see from Equations~\ref{eq:max} and~\ref{eq:min}, the symmetric stair preconditioner results in a smaller condition number as the eigenvalues are of the range:
\begin{equation} \label{eq:finalVals}
    \begin{split}
        \lambda(\Phi_{sym}^{-1}S) &\in (0,1],\\
        \lambda(\Phi_{add}^{-1}S) &\in \left( 0, \frac98 \right].
    \end{split}
\end{equation}

\section{Numerical Results} \label{sec:results}
In this section we present a numerical evaluation of our symmetric stair preconditioner and compare it to the additive stair preconditioner and other parallel preconditioners from the literature on representative trajectory optimization tasks. We evaluate the spread of the eigenvalues both in absolute terms and with regard to the the relative condition number after preconditioning, as well as the resulting number of iterations of PCG needed for convergence. We use the canonical pendulum and cart pole swing up problems as well as a problem to compute a motion across its workspace for a Kuka LBR IIWA-14 manipulator. Source code accompanying this evaluation, including all hyperparameter values used in these experiments, can be found at {\small \href{https://github.com/a2r-lab/SymStair}{\color{blue} \texttt{github.com/a2r-lab/SymStair}}}.

\iffigs
\begin{figure}[!b]
   \centering
   \vspace{-10pt}
   \includegraphics[width=0.99\columnwidth]{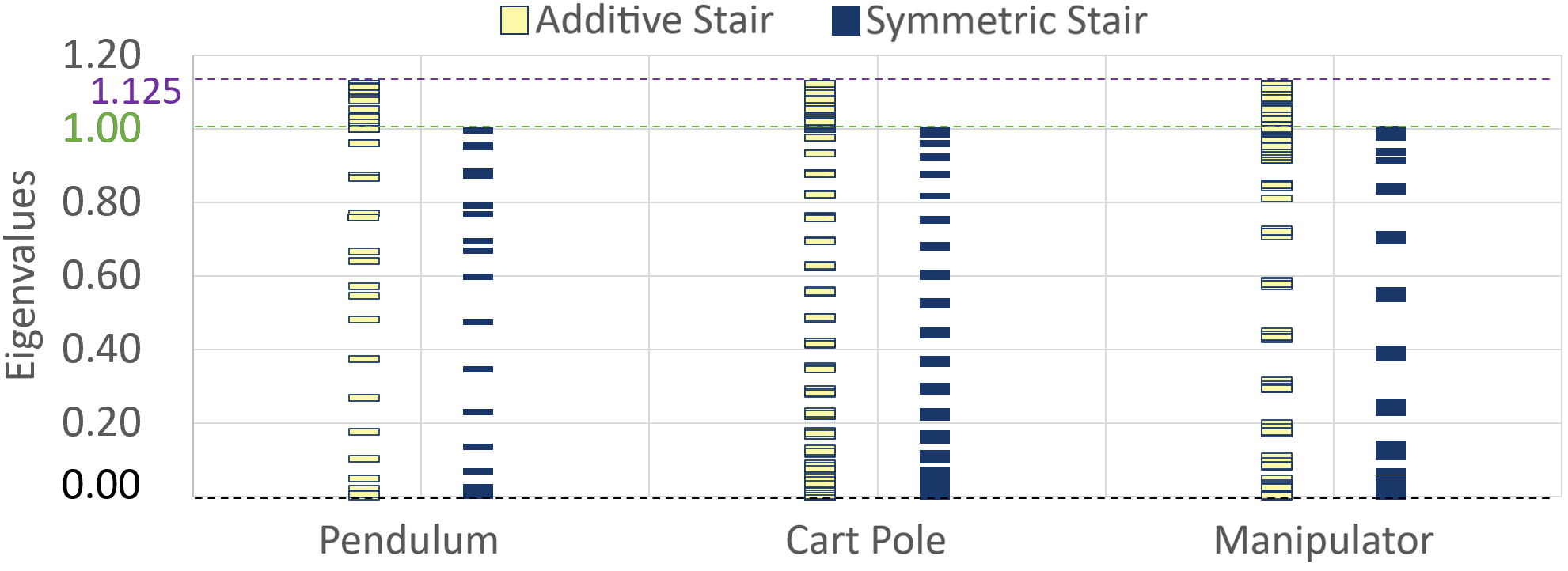}
   \vspace{-20pt}
\caption{Distribution of the Eigenvalues of the additive and symmetric stair preconditioners matching the theoretical results in Equation~\ref{eq:finalVals}.}
\label{fig:eigs}
\end{figure}
\fi

\hyperref[fig:condition]{Figure~\ref{fig:eigs}} shows the distribution of the eigenvalues of the additive and symmetric stair preconditioners along with the line for 0.0 highlighted in black, 1.0 in green, and $\frac{9}{8} = 1.125$ in purple. We see that across all problems our numerical results match the theoretical results in Equation~\ref{eq:finalVals}. 
%

This theoretical result directly leads to a better relative condition number for our preconditioner as compared to alternative parallel preconditioners from the literature. As shown in \hyperref[fig:condition]{Figure~\ref{fig:condition}}, when normalizing the condition number to the Jacobi preconditioner equaling 1 for each system, to enable comparison across the three problems which have large differences in the absolute condition number, we see that the symmetric stair preconditioner is the most performant. In fact, it is not only able to reduce the relative condition number by 76-89\% as compared to the standard Jacobi preconditioner, but also outperforms the best parallel preconditioner from the literature, the additive stair preconditioner, by 33-34\%.

This improvement in numerical conditioning also drastically reduces the number of PCG iterations needed for the problem to converge, enabling faster linear system solves. This is crucial for real-time trajectory optimization as each trajectory optimization solve requires the solution of many linear systems. In particular, as shown in \hyperref[fig:pcg]{Figure~\ref{fig:pcg}}, when solving the first KKT system for our three target trajectory optimization problems, PCG using the symmetric stair preconditioner requires 51-68\% fewer iterations than the Jacobi preconditioner, and 17-25\% fewer iterations than the next best parallel preconditioner to converge to a solution under the same exit tolerance.

\begin{table*}[!b]
\begin{equation} \label{eq:lz}
\scalebox{\scaleForAppendixMath}{$
\Psi_l^{-1}P_l =
    \displaystyle \left(\begin{matrix}0 & - D_{1}^{-1} O_{1} & 0 & 0 & 0 & 0\\0 & D_{2}^{-1} O_{2} D_{3}^{-1} O_{2}^{T} + D_{2}^{-1} O_{1}^{T} D_{1}^{-1} O_{1} & 0 & D_{2}^{-1} O_{2} D_{3}^{-1} O_{3} & 0 & 0\\0 & - D_{3}^{-1} O_{2}^{T} & 0 & - D_{3}^{-1} O_{3} & 0 & 0\\0 & D_{4}^{-1} O_{3}^{T} D_{3}^{-1} O_{2}^{T} & 0 & D_{4}^{-1} O_{4} D_{5}^{-1} O_{4}^{T} + D_{4}^{-1} O_{3}^{T} D_{3}^{-1} O_{3} & 0 & D_{4}^{-1} O_{4} D_{5}^{-1} O_{5}\\0 & 0 & 0 & - D_{5}^{-1} O_{4}^{T} & 0 & - D_{5}^{-1} O_{5}\\0 & 0 & 0 & D_{6}^{-1} O_{5}^{T} D_{5}^{-1} O_{4}^{T} & 0 & D_{6}^{-1} O_{5}^{T} D_{5}^{-1} O_{5}\end{matrix}\right)
$}
\end{equation}
\begin{equation} \label{eq:rz}
\scalebox{\scaleForAppendixMath}{$
    \Psi_r^{-1}P_r =
    \displaystyle \left(\begin{matrix}D_{1}^{-1} O_{1} D_{2}^{-1} O_{1}^{T} & 0 & D_{1}^{-1} O_{1} D_{2}^{-1} O_{2} & 0 & 0 & 0\\- D_{2}^{-1} O_{1}^{T} & 0 & - D_{2}^{-1} O_{2} & 0 & 0 & 0\\D_{3}^{-1} O_{2}^{T} D_{2}^{-1} O_{1}^{T} & 0 & D_{3}^{-1} O_{3} D_{4}^{-1} O_{3}^{T} + D_{3}^{-1} O_{2}^{T} D_{2}^{-1} O_{2} & 0 & D_{3}^{-1} O_{3} D_{4}^{-1} O_{4} & 0\\0 & 0 & - D_{4}^{-1} O_{3}^{T} & 0 & - D_{4}^{-1} O_{4} & 0\\0 & 0 & D_{5}^{-1} O_{4}^{T} D_{4}^{-1} O_{3}^{T} & 0 & D_{5}^{-1} O_{5} D_{6}^{-1} O_{5}^{T} + D_{5}^{-1} O_{4}^{T} D_{4}^{-1} O_{4} & 0\\0 & 0 & 0 & 0 & - D_{6}^{-1} O_{5}^{T} & 0\end{matrix}\right)
$}
\end{equation}

\begin{equation} \label{eq:(l-d)p}
\scalebox{\scaleForAppendixMath}{$
(\Psi_l^{-1} - D^{-1})P_l = \displaystyle \left(\begin{matrix}0 & 0 & 0 & 0 & 0 & 0\\0 & D_{2}^{-1} O_{2} D_{3}^{-1} O_{2}^{T} + D_{2}^{-1} O_{1}^{T} D_{1}^{-1} O_{1} & 0 & D_{2}^{-1} O_{2} D_{3}^{-1} O_{3} & 0 & 0\\0 & 0 & 0 & 0 & 0 & 0\\0 & D_{4}^{-1} O_{3}^{T} D_{3}^{-1} O_{2}^{T} & 0 & D_{4}^{-1} O_{4} D_{5}^{-1} O_{4}^{T} + D_{4}^{-1} O_{3}^{T} D_{3}^{-1} O_{3} & 0 & D_{4}^{-1} O_{4} D_{5}^{-1} O_{5}\\0 & 0 & 0 & 0 & 0 & 0\\0 & 0 & 0 & D_{6}^{-1} O_{5}^{T} D_{5}^{-1} O_{4}^{T} & 0 & D_{6}^{-1} O_{5}^{T} D_{5}^{-1} O_{5}\end{matrix}\right)
$}
\end{equation}

\begin{equation} \label{eq:dpl}
\scalebox{\scaleForAppendixMath}{$
    D^{-1}P_l =  \displaystyle \left(\begin{matrix}0 & - D_{1}^{-1} O_{1} & 0 & 0 & 0 & 0\\0 & 0 & 0 & 0 & 0 & 0\\0 & - D_{3}^{-1} O_{2}^{T} & 0 & - D_{3}^{-1} O_{3} & 0 & 0\\0 & 0 & 0 & 0 & 0 & 0\\0 & 0 & 0 & - D_{5}^{-1} O_{4}^{T} & 0 & - D_{5}^{-1} O_{5}\\0 & 0 & 0 & 0 & 0 & 0\end{matrix}\right) \quad \quad
    (\Psi_l^{-1} - D^{-1})\Psi_l =  \displaystyle \left(\begin{matrix}0 & 0 & 0 & 0 & 0 & 0\\- D_{2}^{-1} O_{1}^{T} & 0 & - D_{2}^{-1} O_{2} & 0 & 0 & 0\\0 & 0 & 0 & 0 & 0 & 0\\0 & 0 & - D_{4}^{-1} O_{3}^{T} & 0 & - D_{4}^{-1} O_{4} & 0\\0 & 0 & 0 & 0 & 0 & 0\\0 & 0 & 0 & 0 & - D_{6}^{-1} O_{5}^{T} & 0\end{matrix}\right)
$}
\end{equation}
\end{table*}

\iffigs
\begin{figure}[!t]
   \centering
   \vspace{5pt}
   \includegraphics[width=\columnwidth]{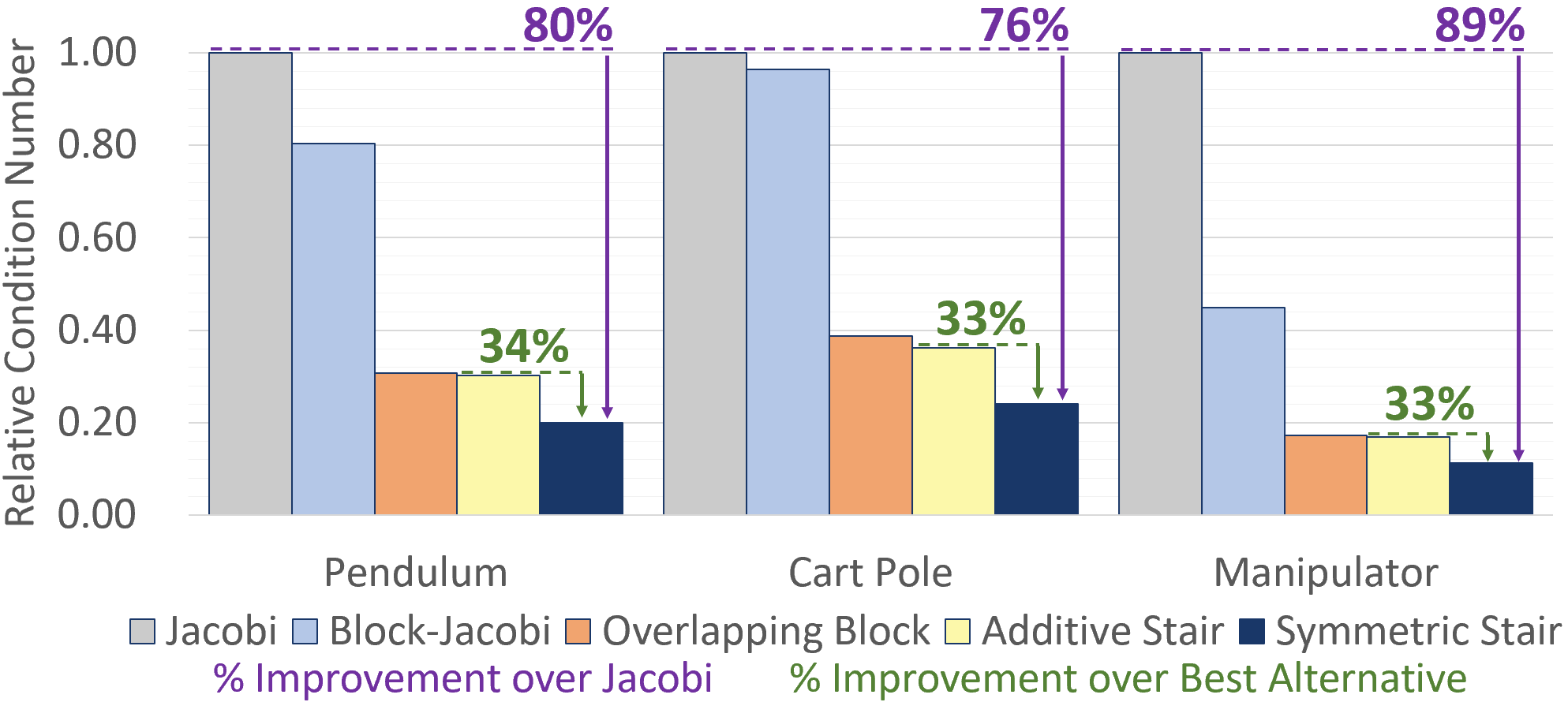}
   \vspace{-20pt}
\caption{The relative condition number resulting from different preconditioners, normalized to the value of the Jacobi preconditioner, showing the improved performance of the symmetric stair preconditioner.}
\vspace{-10pt}
\label{fig:condition}
\end{figure}
\fi

\iffigs
\begin{figure}[!t]
   \centering
   \vspace{5pt}
   \includegraphics[width=\columnwidth,trim={0 0 0 4pt}]{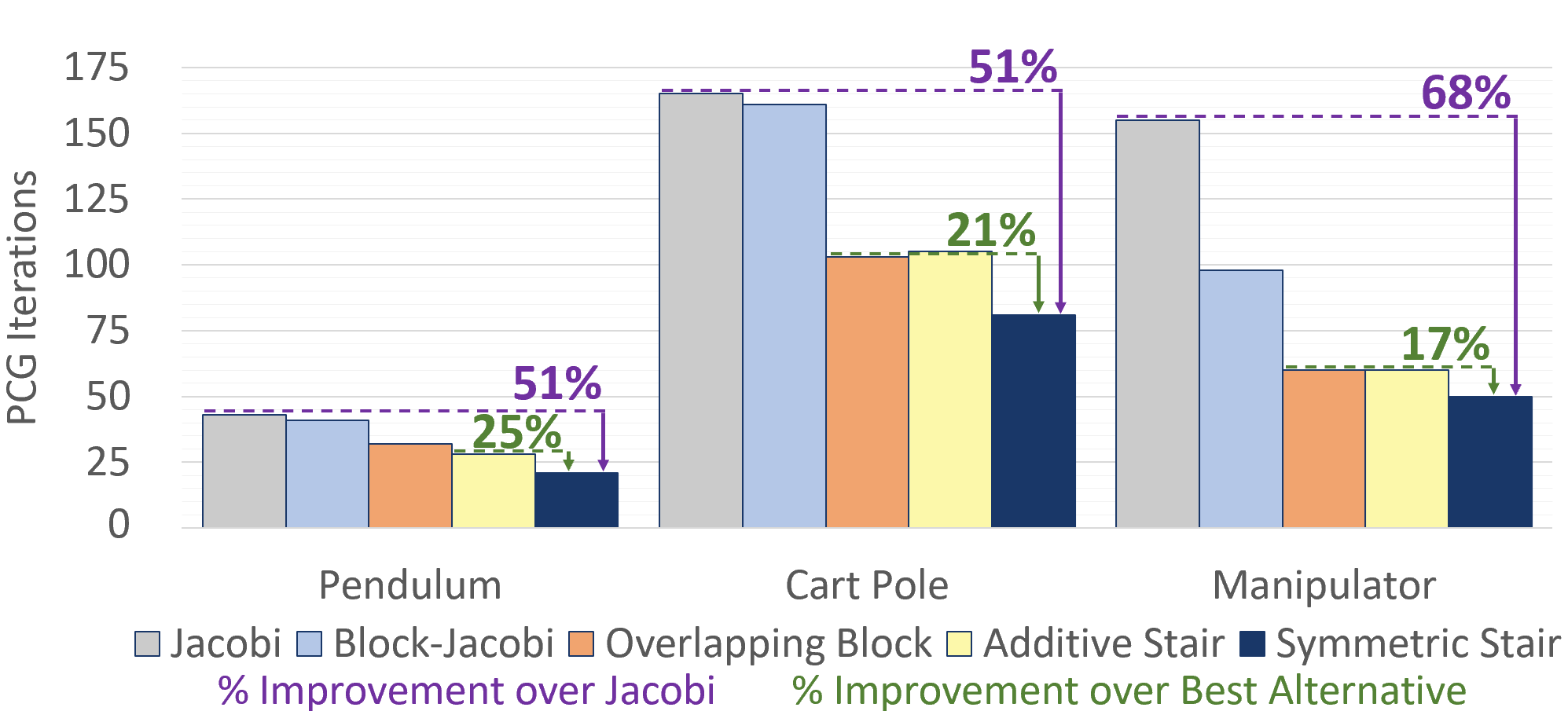}
   \vspace{-20pt}
\caption{The number of PCG iterations required for convergence to the same exit tolerance across different preconditioners and problems, again showing the improved performance of the symmetric stair preconditioner.}
\vspace{-10pt}
\label{fig:pcg}
\end{figure}
\fi

\section{Conclusion and Future Work} \label{sec:conclusion}
In this work we present a new parallel-friendly symmetric stair preconditioner. Through both proofs and numerical experiments, we show that our preconditioner has advantageous theoretical and practical properties that enable fast iterative linear system solves. Across three benchmark tasks, our preconditoner provides up to a 34\% reduction in condition number and 25\% reduction in PCG iterations. 

In future work we hope to explore the theoretical and practical performance tradeoffs of higher order polynomial preconditioners on heterogeneous hardware, as well as extend our analysis to provide bounds on the number of PCG iterations under various exit tolerances and preconditioners.

\section{Appendix} \label{sec:Appendix}
Some useful results revealing the sparsity patterns for the $6\times 6$ block cases are listed in Equations~\ref{eq:lz}-~\ref{eq:dpl}.

\clearpage

\bibliographystyle{bib/IEEEtran_new}
\bibliography{bib/refs.bib}

\end{document}